\documentclass[11pt]{article}
\usepackage{amsmath,amsfonts,amssymb,amsthm} 
\usepackage[explicit]{titlesec}
\usepackage{graphicx}
\usepackage{enumitem}
\setlist{nolistsep}
\usepackage{stmaryrd}
\usepackage{mathpple}
\usepackage{enumerate}
\usepackage{textcomp}
\usepackage[utf8]{inputenc}
\usepackage{mathtools}
\usepackage{epsfig}
\usepackage{xcolor}
\usepackage{array}
\usepackage{calc}
\usepackage{fancyhdr}
\usepackage[colorlinks=false]{hyperref}
\usepackage[noabbrev]{cleveref}
\usepackage{dsfont}
\usepackage{wrapfig}
\usepackage{fancyhdr}
\usepackage{float}
\usepackage[font={small}]{caption}
\usepackage{eurosym}
\usepackage{csquotes}
\usepackage[left=3cm,right=3cm,top=2cm,bottom=2.5cm]{geometry}

\DeclareMathAlphabet{\mathcal}{OMS}{cmsy}{m}{n}
\newcommand{\R}{\mathbb{R}}
\newcommand{\IP}{\mathbb{P}}

\newcommand{\N}{\mathbb{N}}

\newcommand{\E}{\mathbb{E}}

\newcommand{\IL}{\mathbb{L}}

\newcommand{\cB}{\mathcal{B}}
\newcommand{\cF}{\mathcal{F}}

\newcommand{\cU}{\mathcal{U}}
\newcommand{\cR}{\mathcal{R}}
\newcommand{\cN}{\mathcal{N}}
\newcommand{\cM}{\mathcal{M}}

\newcommand{\cW}{\mathcal{W}}

\newcommand{\1}{\mathds{1}}

\newcommand{\drm}{\mathrm{d}}

\newcommand{\eps}{\varepsilon}
\newcommand{\cMc}{\langle M \rangle}

\newcommand{\floor}[1]{\lfloor #1 \rfloor}

\newcommand{\comb}[2]{\Big({\scriptstyle\begin{matrix}#1 \\ #2\end{matrix}}\Big)}
\newcommand{\blue}[1]{{\color{black}#1}}

\font\calcal=cmsy10 scaled\magstep1
\def\build#1_#2^#3{\mathrel{\mathop{\kern 0pt#1}\limits_{#2}^{#3}}}
\def\liml{\build{\longrightarrow}_{}^{{\mbox{\calcal L}}}}

\numberwithin{equation}{section}

\setitemize[1]{label=--,partopsep=\parskip,topsep=-\parskip,itemsep=0pt,parsep=0pt}
\titleformat{\chapter}{\fontfamily{phv}\selectfont\LARGE\bfseries}{\thesection}{1em}{\fontfamily{phv}\selectfont\LARGE\bfseries #1}
\titleformat{\section}{\centering\fontfamily{phv}\selectfont\bfseries}{\thesection}{1em}{\centering\fontfamily{phv}\selectfont\bfseries #1}
\titleformat{\subsection}{\fontfamily{phv}\selectfont\small\bfseries\centering}{\thesubsection}{1em}{\fontfamily{phv}\selectfont\small\bfseries #1}
\titleformat{\subsubsection}{\fontfamily{phv}\selectfont\footnotesize\bfseries\centering}{\thesubsubsection}{1em}{\fontfamily{phv}\footnotesize\selectfont\bfseries #1}

\renewcommand*\thesection{\arabic{section}}

\makeatletter
\renewcommand{\@seccntformat}[1]{\llap{{\csname the#1\endcsname}\hspace{1em}}}   

\newtheoremstyle{thmit}{10pt}{10pt}
{\normalfont\itshape}
{}
{\small\bf\fontfamily{phv}\selectfont}
{\;}
{0.25em}
{\small\fontfamily{phv}\selectfont\thmname{#1}\nobreakspace\footnotesize\thmnumber{#2}.
\thmnote{\nobreakspace\the\thm@notefont\fontfamily{phv}\selectfont\footnotesize\bfseries\--\nobreakspace#3}}

\newtheoremstyle{rmq}
{5pt}
{5pt}
{\normalfont}
{}
{\small\bf\fontfamily{phv}\selectfont}
{\;}
{0.25em}
{\small\fontfamily{phv}\selectfont\thmname{#1}\nobreakspace\footnotesize\thmnumber{\@ifnotempty{#1}{}\@upn{#2}}
\thmnote{\nobreakspace\the\thm@notefont\fontfamily{phv}\selectfont\footnotesize\bfseries--\nobreakspace#3.}}
\makeatother

\newcounter{count}
\numberwithin{count}{section}
\newcounter{alpha}

\theoremstyle{thmit}
\newtheorem{theorem}[count]{Theorem\small}

\newtheorem{lemma}[count]{Lemma\small}
\newtheorem{corollary}[count]{Corollary\small}
\newtheorem{remark}[count]{Remark\small}

\theoremstyle{rmq}

\renewenvironment{proof}[2]
{\paragraph{\fontfamily{phv}\selectfont\bfseries\small Proof of #1 \footnotesize #2.}}%
{\begin{flushright}
\qed
\end{flushright}}

\begin{document}
\pagestyle{fancy}

\title{\vspace{-3ex}
\fontfamily{phv}\selectfont\bfseries\Large
\blue{New insights on the reinforced Elephant Random Walk using a martingale approach}}
\author{\fontfamily{phv}\selectfont\bfseries
Lucile Laulin}
\date{}
\AtEndDocument{\bigskip{\footnotesize%
\textsc{Université de Bordeaux, Institut de Mathématiques de Bordeaux,
UMR 5251, 351 Cours de la Libération, 33405 Talence cedex, France.} \par  
\textit{E-mail adress :} \href{mailto:lucile.laulin@math.u-bordeaux.fr}{\texttt{lucile.laulin@math.u-bordeaux.fr}} \par
}}

\maketitle

\centerline{
\begin{minipage}[c]{0.7\textwidth}
{\small\section*{Abstract}\vspace{-1ex}
This paper is devoted to the asymptotic analysis of the reinforced elephant random walk (RERW) using a martingale approach. In the diffusive and critical regimes, we establish the almost sure convergence, the law of iterated logarithm and the quadratic strong law for the RERW. The distributional convergences of the RERW to some Gaussian processes are also provided. In the superdiffusive regime, we prove the distributional convergence as well as the mean square convergence of the RERW. All our analysis relies on asymptotic results for multi-dimensional martingales with matrix normalization.
}\medskip\\
\small
{\bf MSC:} primary 60G50; secondary 60G42; 60F17\medskip\\
{\bf Keywords : }Elephant random walk; Reinforced random walk; Multi-dimensional martingales; Almost sure convergence; Asymptotic normality; Distributional convergence
\end{minipage}
}

\setlength{\parindent}{0pt}

\section{Introduction}\thispagestyle{empty}
Reinforced random walks have generated much interest in the recent years with the focus being mainly on graphs, edge or vertex reinforced random walk, see for example \cite{Kozam2013} or \cite{Pemantle2007} for a comprehensive and extensive overview on the subject, as well as the recent contribution \cite{baur2019,Bertoin2020}.

\blue{In this paper, we investigate a special case of reinforced random walk in connection with the Elephant Random Walk (ERW)}, introduced by Schütz and Trimper \cite{Schutz2004} in the early 2000s. At first, the ERW was used in order to see how long-range memory affects the random walk and induces a crossover from a diffusive to superdiffusive behavior. It was referred to as the ERW in allusion to the traditional saying that elephants can always remember anywhere they have been.  
The elephant starts at the origin at time zero, $S_0 = 0$. At time $n = 1$,  the elephant moves in one to the right with probability $q$ and to the left with probability $1-q$ for some $q$ in $[0,1]$.  
Afterwards, at time $n+1$, the elephant chooses uniformly at random an integer $k$ among the previous times $1,\ldots,n$. Then, it moves exactly in 
the same direction as that of time $k$ with probability $p$ or the oppositve direction with the probability $1-p$, where the parameter $p$ 
stands for the memory parameter of the ERW. The position of the elephant at time $n+1$ is given by
\begin{equation} 
\label{POS-bERW}
S_{n+1} = S_n + X_{n+1}
\end{equation}
where $X_{n+1}$ is the $(n+1)$-th increment of the random walk.
The ERW shows three differents regimes depending on the location of its memory parameter $p$ with respect to the critical value
$p_c=3/4$.

A wide literature is now available on the ERW in dimension $d=1$. A strong law of large numbers and a central limit theorem for the position $S_n$, properly normalized, were established in the diffusive regime $p< 3/4$ and the critical regime $p=3/4$, see \cite{Baur2016}, \cite{Coletti2017}, \cite{ColettiN2017}, \cite{Schutz2004} and the more recent contributions \cite{BercuHG2019}, \cite{Coletti2019}, \cite{Fan2020}, \cite{Gonzales2020}, \cite{Takei2020}, \cite{Vazquez2019}. 
The superdiffusive regime $p>3/4$ turns out to be harder to deal with.
Bercu \cite{Bercu2018} proved that the limit of the position of the ERW is not Gaussian and Kubota and Takei \cite{Kubota2019} showed that 
the fluctuation of the ERW around its limit in the superdiffusive regime is Gaussian. Finally, Bercu and Laulin in \cite{BercuLaulin2019}
extended all the results of \cite{Bercu2018} to the multi-dimensional ERW (MERW) where $d \geq 1$ and to its center of mass \cite{LaulinBercu2020}. Moreover, functional central limit theorems were also provided via a connection to Pólya-type urns, see \cite{Baur2016} for the ERW, \cite{baur2019} for a particular class of random walks with reinforced memory such as the ERW and the Shark Random Swim \cite{Businger2018},  and more recently \cite{Bertenghi2020} for the MERW.
\medskip\\
The main subjet of this paper is to study the asymptotical behavior of the reinforced Elephant Random Walk (RERW).
As it was done in \cite{Bercu2018}, we can write the $(n+1)$-th increment $X_{n+1}$ under the form 
\begin{equation}
\label{X-a-b}
X_{n+1}=\alpha_{n+1} X_{\beta_{n+1}}.
\end{equation}
In the case of the ERW, we had $\alpha_{n+1}\sim\cR(p)$ and $\beta_{n+1}\sim \cU\{1,\ldots,n\}$. The major change for the RERW is that the distribution of $\beta_n$ is no longer uniform.
\smallskip\\
\blue{Very recently, Baur \cite{baur2019} studied the asymptotic behavior of the RERW using a Polya-type urns approach. He established interesting functional limit theorems thanks to the seminal work of Janson \cite{Janson2004}. Our strategy is totally different as it relies on a martingale approach. On the one hand, we prove new almost sure convergence results such as strong laws of large numbers, laws of iterated logarithms, as well as quadratic strong laws. On the other hand, we give an alternative method to obtain the functional limit theorems without making use of the results from \cite{Janson2004}. The martingale approach we propose fulfills these two objectives. The main strength of our approach is that calculations are totally sefl-contained and rather easy to follow. It should also be noted that using the martingale theory is sufficient on its own to obtain all the results presented in this paper. Moreover, we strongly believe that this could be used to study several variations of the ERW with reinforced memory or more generally reinforced random walks.} 
 
\smallskip
This paper is organized as follows. The model of reinforced memory is presented in Section \ref{S-model} while the main results are given in Section \ref{S-main-results}. We first investigate the diffusive regime and we establish the almost sure convergence, the law of iterated logarithm and the quadratic strong law for the RERW. The functional central limit theorem is also provided. Next, we prove similar results in the critical regime. Finally, we establish a strong limit theorem in the superdiffusive regime. Our martingale approach is described in Section \ref{S-martingale-approach}. Finally, all technical proofs are
postponed to Sections \ref{S-AS-results}--\ref{S-FCLT-results}.

\section{The reinforced elephant random walk}
\label{S-model}
We assume in all the sequel that the memory parameter $p\neq 1/2$ since the particular case $p=1/2$ reduces to the standard random walk.
\blue{Let $\cF_n=\sigma(X_1,\ldots,X_n,\beta_1,\ldots,\beta_n)$ be the natural $\sigma$-algebra up to time $n$ and denote by $\rho_{n}(k)$ the weight of the instant $k$ after $n$ steps.} The ERW is associated with the special case where $\rho_{n}(k)=1$ if $k\leq n$ and $0$ elsewise. Adding a reinforcement of weight $c$, where $c$ is a non-negative real number, implies that the weight $\rho_{n}(k)$ of instant $k$ is modified as follows
\begin{equation*}
\rho_{n}(k)=\left\{\begin{array}{ll}
0 & \text{if $k\geq n+1$}, \\
1 & \text{if $k=n$}, \\
\rho_{n-1}(k)+ c \1_{\beta_{n}=k} & \text{if $1\leq k<n$}.
\end{array}\right.
\end{equation*}
Consequently, it follows from the very definition of $\rho_n(k)$ that the conditional distribution of $\beta_{n+1}$ is given by, for $1\leq k\leq n$,
\begin{equation*}
\IP(\beta_{n+1}=k|\blue{\cF_n})=\frac{\rho_{n}(k)}{\sum_{j=1}^n \rho_{n}(j)}=\frac{\rho_{n}(k)}{(c+1)n-c}.
\end{equation*}
The parameter $c$ represents the intensity of the reinforcement. The reader can notice that the case $c=0$ corresponds to the traditionnal ERW, and that in \blue{this case the distribution of $\beta_{n+1}$ is only dependant of the time $n$}.
Hereafter, let $a=2p-1$, such that $-1 \leq a \leq 1$. We have by the definition of $X_n$,
\begin{equation*}
\E[X_{n+1}|\cF_n]=\E[\alpha_{n+1}]\E[X_{\beta_{n+1}}|\cF_n]\nonumber  =a\E\Big[\sum_{k=1}^n X_k\1_{\beta_{n+1}=k}|\cF_n\Big]\nonumber =\frac{a}{(c+1)n-c}\sum_{k=1}^n X_k\rho_{n}(k)\nonumber.
\end{equation*}
Then, denote
\begin{equation}
Y_n =\sum_{k=1}^n X_k\rho_{n}(k)
\end{equation}
such that $Y_n = S_n$ when $c=0$, and
\begin{equation}
\label{ESPC-X}
\E[X_{n+1}|\cF_n] =\frac{a}{(c+1)n-c}Y_n.
\end{equation}
Hence, we immediately get
\begin{equation}
\label{ESPC-S}
\E[S_{n+1}|\cF_n] = S_n + \E[X_{n+1}|\cF_n] = S_n + \frac{a}{(c+1)n-c}Y_n.
\end{equation}
Hereafter, notice that
\begin{align}
\label{YN-YN1-XB}
Y_{n+1}  =\sum_{k=1}^{n+1}X_k\rho_{n+1}(k) 
 = \sum_{k=1}^{n}X_k\big(\rho_{n}(k)+c \1_{\beta_{n+1}=k}\big) + X_{n+1} = Y_n + (\alpha_{n+1}+c) X_{\beta_{n+1}} 
\end{align}
we obtain
\begin{equation}
\label{ESPC-Y}
\E[Y_{n+1}|\cF_n]  = \Big(1 + \frac{a+c}{(c+1)n-c}\Big)Y_n.
\end{equation}
Finally, for any $n\geq 1$ let
\begin{equation}
\label{DEF-gammaN}
\gamma_n= 1+\frac{a+c}{(c+1)n-c}=\frac{n+a\lambda}{n-c\lambda} 
\quad\text{where}\quad\lambda=\frac{1}{c+1}
\end{equation}
and 
\begin{equation}
\label{DEF-aN}
a_{n}=\prod_{k=1}^{n-1}\gamma_k^{-1} =\frac{\Gamma(n-c\lambda)\Gamma(1+a\lambda)}{\Gamma(n+a\lambda)\Gamma(\lambda)}.
\end{equation}
It follows from standard calculations on the Gamma function that
\begin{equation}
\label{lim-CN}
\lim_{n\to\infty} n^{(a+c)\lambda}a_n =  \frac{\Gamma(1+a\lambda)}{\Gamma(\lambda)}.
\end{equation}
Our strategy for proving asymptotic results for the reinforced elephant random
walk is as follows. On the one hand, the behavior of position $S_n$ is closely
related to the one of the sequences $(M_n)$ and $(N_n)$ defined for all $n\geq 0$ by
\begin{equation}
\label{DEF-M-N}
M_n = a_n Y_n \quad\text{and}\quad N_n = S_n - \frac{a}{a+c}Y_n.
\end{equation}
We immediately get from \eqref{ESPC-Y} and \eqref{DEF-aN} that $(M_n)$ is a locally square-integrable martingale adapted to $\cF_n$. Moreover, we have from \eqref{ESPC-X},\eqref{ESPC-S} and \eqref{ESPC-Y} that
\begin{equation*}
E\Big[S_{n+1}-\frac{a}{a+c}Y_{n+1}|\cF_n\Big]=S_n - \frac{a}{a+c}Y_n
\end{equation*}
which means that $(N_n)$ is also a locally square-integrable martingale adapted to $\cF_n$.
On the other hand, we can rewrite $S_n$ as
\begin{equation}
\label{SN-MG}
S_n = N_n + \frac{a}{a+c}a_n^{-1}M_n
\end{equation}
and equation \eqref{SN-MG} allows us to establish the asymptotic behavior of the RERW via an extensive use of the strong law of large numbers and the functional central limit theorem for multi-dimensional martingales \cite{Chaabane2000}, \cite{Duflo1997}, \cite{Hall1980}, \cite{Touati1991}.

\section{Main results}
\label{S-main-results}
\subsection{The diffusive regime}
Our first result deals with the strong law of large numbers for the RERW in the diffusive regime where $a<(1-c)/2$.
\begin{theorem}
\label{THM-CVPS-DIF}
We have the almost sure convergence
\begin{equation}
\lim_{n\to\infty} \frac{S_n}{n}=0 \quad\text{a.s.}
\end{equation}
\end{theorem}
The almost sure rate of convergence for RERW is as follows.
\begin{theorem}
\label{THM-LFQ-DIF}
We have the quadratic strong law
\begin{equation}
\label{QUAD-DIF}
\lim_{n \rightarrow \infty} \frac{1}{\log n} \sum_{k=1}^n \frac{S_k^2}{k^2}=\frac{2ac+c-1}{2a+c-1}\quad  \text{a.s.}
\end{equation}

\begin{remark}
In addition, we could also obtain an upper-bound for the law of iterated logarithm as it was done for the center of mass of the MERW in \cite{LaulinBercu2020}.
\end{remark}

\end{theorem}
Hereafter, we are interested in the distributional convergence of the RERW, which holds in the Skorokhod space $D([0,\infty[)$ of right-continuous functions with left-hand limits. \blue{The following theorem was first obtained by Baur \cite[Theorem 3.2]{baur2019} in the case of a memory parameter equal to $(p+1)/2$.}
\begin{theorem}
\label{CLFT-RERW-DR}
The following convergence in distribution in $D([0,\infty[)$ holds
\begin{equation}
\label{CVL-cont-DIF}
\Big(\frac{S_{\floor{nt}}}{\sqrt{n}}, \ t\geq 0\Big) \Longrightarrow \big(W_t, \ t\geq 0\big)
\end{equation}
where $\big(W_t, \ t\geq 0\big)$ is a real-valued centered Gaussian process starting from the origin with covariance
\begin{equation}
\label{ESP-W}
\E[W_sW_t] = \frac{a(1-c^2)}{(a+c)(1-2a-c)}s\Big(\frac{t}{s}\Big)^{\lambda(a+c)}+ \frac{c(a+1)}{a+c}s
\end{equation}
for $0<s\leq t$. In particular, we have
\begin{equation}
\label{CVL-DIF}
\frac{S_n}{\sqrt{n}} \liml \cN \left(0, \frac{2ac+c-1}{2a+c-1}\right).
\end{equation}
\end{theorem}
\begin{remark}
When $c=0$ we find again the results from \cite{Baur2016} for the ERW
\begin{equation*}
\Big(\frac{S_{\floor{nt}}}{\sqrt{n}}, \ t\geq 0\Big) \Longrightarrow \big(W_t, \ t\geq 0\big)
\end{equation*}
where $\big(W_t, \ t\geq 0\big)$ is a real-valued mean-zero Gaussian process starting from the origin and
\begin{equation*}
\label{ESP-W}
\E[W_sW_t] = \frac{1}{1-2a}s\Big(\frac{t}{s}\Big)^{a}.
\end{equation*}
In particular, we also obtain the asymptotic normality from \cite{Bercu2018,Coletti2017}
\begin{equation*}
\label{def-Gn}
\frac{S_n}{\sqrt{n}} \liml \cN \left(0, \frac{1}{1-2a}\right).
\end{equation*}
\end{remark}
As it was done in \cite{Bertenghi2020}, we also obtain the asymptotic normality for the center of mass of the RERW defined by $$G_n = \frac{1}{n}\sum_{k=1}^n S_n.$$
\begin{corollary}
\label{coro-Gn}
We have the asymptotic normality
\begin{equation}
\frac{G_n}{\sqrt{n}} \liml \cN \left(0, \frac{2 - c(c+1 +3ca +3a -2a^2)}{3(2+c-a)(1-2a-c)}\right).
\end{equation}
\end{corollary}
\begin{remark}
When $c=0$, we find again the asympotic normality established in \cite{LaulinBercu2020,Bertenghi2020}
\begin{equation*}
\frac{G_n}{\sqrt{n}} \liml \cN \left(0, \frac{2}{3(1-2a)(2-a)}\right).
\end{equation*}
\end{remark}
\begin{center}
\begin{minipage}[c]{0.48\linewidth}
\includegraphics[width=8cm]{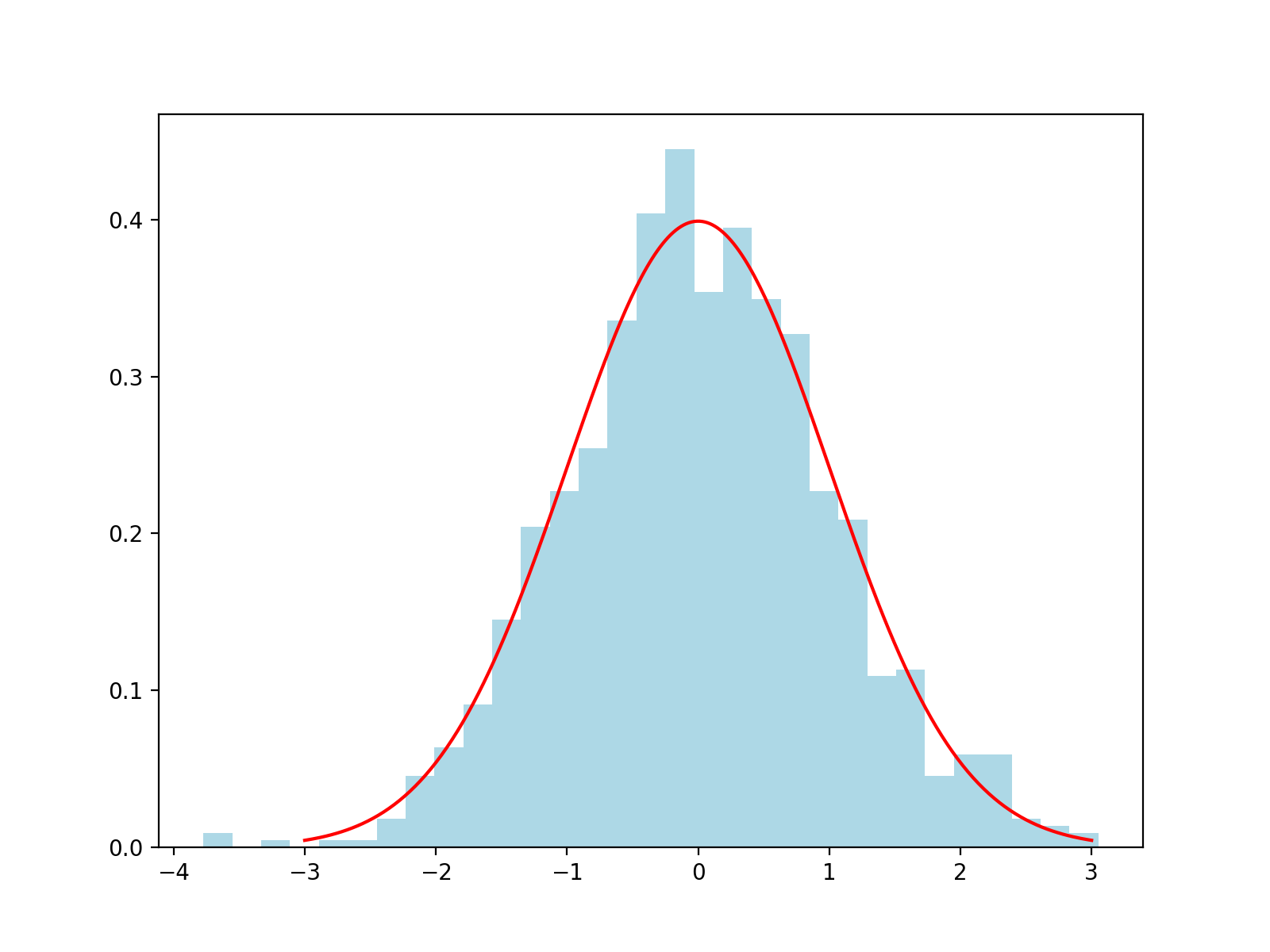}
\captionof{figure}{Asymptotic normality for the RERW in the diffusive regime, when $p=0.35$ and $c=1$.}
\label{picture-TLC-RERW1}
\end{minipage}\hfill
\end{center}
\subsection{The critical regime}
Hereafter, we investigate the critical regime where $a=(1-c)/2$.
\begin{theorem}
\label{THM-CVPS-CR}
We have the almost sure convergence
\begin{equation}
\lim_{n\to\infty} \frac{S_n}{\sqrt{n}\log n}=0\quad\text{a.s.}
\end{equation}
\end{theorem}
The almost sure rates of convergence for the RERW are as follows.
\begin{theorem}
\label{THM-LFQ-LIL-CR}We have the quadratic strong law
\begin{equation}
\label{QUAD-CR}
\lim_{n \rightarrow \infty} \frac{1}{\log \log n} \sum_{k=1}^n \frac{S_k^2}{(k\log k)^2}=\frac{(c-1)^2}{c+1}\quad  \text{a.s.}
\end{equation}
In addition, we also have the law of iterated logarithm
\begin{equation}
\label{LIL-CR}
\limsup_{n \rightarrow \infty} \frac{S_n^2}{2 n \log n \log \log \log n}= \frac{(c-1)^2}{c+1}  \quad \text{a.s.}
\end{equation}
\end{theorem}
Once again, our next result concerns the functional convergence in distribution for the RERW. \blue{The following theorem was also first obtained by Baur \cite[Theorem 3.2]{baur2019}.}
\begin{theorem}
\label{CLFT-RERW-CR}
The following convergence in distribution in $D([0,\infty[)$ holds
\begin{equation}
\label{CVL-cont-CR}
\Big(\frac{S_{\floor{n^t}}}{\sqrt{n^t\log n}}, \ t\geq 0\Big) \Longrightarrow \sqrt{\frac{(c-1)^2}{(c+1)}}\big(B_t, \ t\geq 0\big)
\end{equation}
where $(B_t, \ t\geq0)$ is a one-dimensional standard Brownian motion. In particular, we have
\begin{equation}
\label{CVL-CR}
\frac{S_n}{\sqrt{n\log n}} \liml \cN \left(0, \frac{(c-1)^2}{c+1} \right).
\end{equation}
\end{theorem}
\begin{remark}
When $c=0$, we find again the results from \cite{Baur2016} for the ERW
\begin{equation*}
\Big(\frac{S_{\floor{nt}}}{\sqrt{n^t\log n}}, \ t\geq 0\Big) \Longrightarrow \big(B_t, \ t\geq 0\big)
\end{equation*}
where $(B_t,\ t\geq0)$ is a one-dimensional standard Brownian motion.
In particular, we find once again the asymptotic normality from \cite{Baur2016,Bercu2018,Coletti2017}
\begin{equation*}
\frac{S_n}{\sqrt{n\log n}} \liml \cN \left(0, 1\right).
\end{equation*}
\end{remark}
\begin{center}
\begin{minipage}[c]{0.48\linewidth}
\includegraphics[width=8cm]{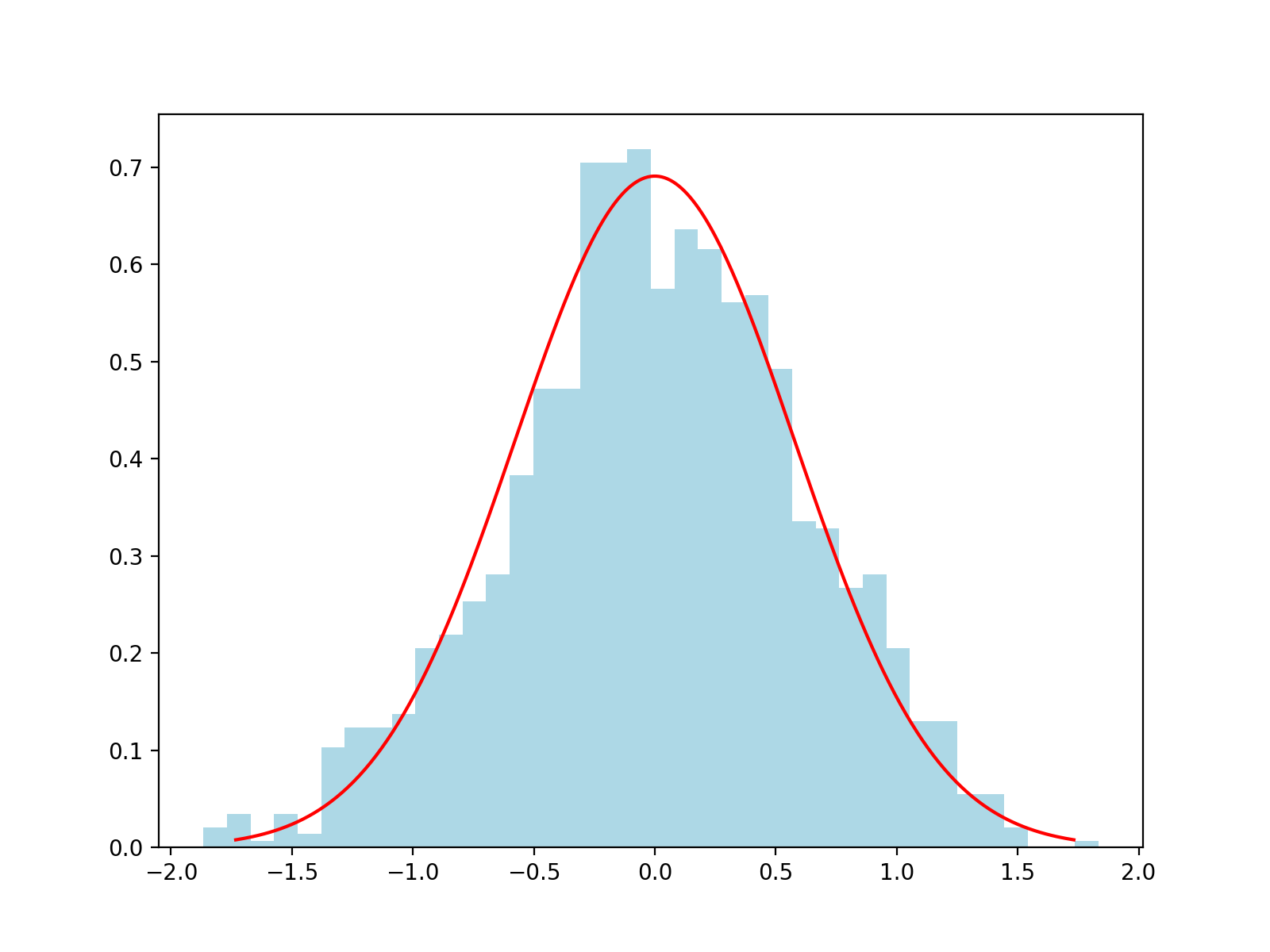}
\captionof{figure}{Asymptotic normality for the RERW in the critical regime, when $c=2$ (ie. $p=0.25$).}
\label{picture-TLC-RERW2}
\end{minipage}
\end{center}
\subsection{The superdiffusive regime}
Finally, we focus our attention on the superdiffusive regime where $a>(1-c)/2$. The reader can notice that it is the only type of behavior for the RERW that still holds when $c>3$ since $a\geq -1$. \blue{The following convergence in $D([0,\infty[)$ can also be found in \cite[Theorem 3.2]{baur2019}. The almost sure and mean-square convergences are new.}
\begin{theorem}
\label{THM-CVPS-SDIF}
We have the following distributional convergence in $D([0,\infty[)$
\begin{equation}
\label{CVPS-cont-SDIF}
\Big(\frac{S_{\floor{nt}}}{n^{\lambda(c+a)}}, \ t\geq 0\Big) \Longrightarrow (\Lambda_t, \ t\geq 0)
\end{equation}
where the limiting $\Lambda_t = t^{\lambda(c+a)}L_c$, $L_c$ being some non-denegerate random variable. In particular, we have
\begin{equation}
\label{CVPS-SDIF}
\lim_{n\to\infty} \frac{S_n}{n^{(a+c)\lambda}}=L_c \quad\text{a.s.}
\end{equation}
Moreover, we also have the mean square convergence
\begin{equation}
\label{CVL2-SDIF}
\lim_{n\to\infty} \E\Big[\Big|\frac{S_n}{n^{(a+c)\lambda}}-L_c\Big|^2\Big]=0.
\end{equation}
\end{theorem}
\begin{theorem}
\label{THM-L2-SDIF}
The expected value of $L_c$ is 
\begin{equation}
\label{ESP-SDIF}
\E[L_c]= \frac{a(2q-1)\Gamma(\lambda)}{(a+c)\Gamma(1+a\lambda)}
\end{equation} 
while its variance is given by
\begin{equation}
\label{L2-SDIF}
\E\big[L_c^2\big]= \frac{a^2(1+2ac+c^2)\Gamma(\lambda)}{(a+c)^2\lambda(2a+c-1)\Gamma((2a+c)\lambda)}.
\end{equation}
\end{theorem}
\begin{remark}
When $c=0$, we find once again the moments of $L$ established in \cite{Bercu2018}
\begin{equation*}
\E[L]= \frac{2q-1}{\Gamma(a+1)}
\quad\text{and}\quad
\E\big[L^2\big]= \frac{1}{(2a-1)\Gamma(2a)}.
\end{equation*}
\end{remark}

\section{A two-dimensional martingale approach}
\label{S-martingale-approach}
In order to investigate the asymptotic behavior of $(S_n)$, we introduce the two-dimensional martingale $(\cM_n)$ defined by
\begin{equation}
\label{DEF-MMN}
\cM_n = \begin{pmatrix} N_n \\ M_n\end{pmatrix}
\end{equation}
where $(M_n)$ and $(N_n)$ are the two locally square-integrable martingales introduced in \eqref{DEF-M-N}. As for the center of mass of the ERW \cite{LaulinBercu2020}, the main difficulty we face is that the predictable quadratic variation of $(M_n)$ and $(N_n)$ increase
to infinity with two different speeds. A matrix normalization will again be necessary to establish the asymptotic behavior of the RERW. \blue{We will study $(\cM_n)$, instead of $(M_n)$ or $(N_n)$.}
\medskip

Let $\eps_{n+1} = Y_{n+1}-\gamma_nY_n$ and $\xi_n =(\alpha_{n}-{a})X_{\beta_{n}}$.
We have from equations \eqref{YN-YN1-XB}, \eqref{DEF-aN} and \eqref{DEF-M-N}
\begin{align}
\label{quad-Mn}
\Delta\cM_{n+1} & =\cM_{n+1}-\cM_n\nonumber \\
			& = \begin{pmatrix}
				S_{n+1}-S_n - \frac{a}{a+c}\big(Y_{n+1}-Y_n\big) \\
					{a_{n+1}Y_{n+1}-a_nY_n}
			\end{pmatrix}\nonumber\\
			& = \begin{pmatrix}\alpha_{n+1}X_{\beta_{n+1}} - \frac{a}{a+c}(\alpha_{n+1}+c)X_{\beta_{n+1}} \\
			{a_{n+1}\eps_{n+1}}
			\end{pmatrix}\nonumber \\
			& = a_{n+1}\eps_{n+1}\comb{0}{1} +\frac{c}{a+c}\xi_{n+1}\comb{1}{0}.
\end{align}

We also find from \eqref{YN-YN1-XB} that
\begin{align}
\label{ESPC-eps2}
\E[\eps_{n+1}^2|\cF_n]&= \E[Y_{n+1}^2|\cF_n]-\gamma_n^2Y_n^2 \nonumber\\
& = Y_n^2 + 2 (\gamma_n-1)Y_n^2 + 1 +2ac+c^2 - \gamma_n^2 Y_n^2 \nonumber\\
& = 1 +2ac+c^2 - (\gamma_n-1)^2Y_n^2.
\end{align}
In addtion, we obtain once again from \eqref{YN-YN1-XB} that
\begin{equation}
\label{ESPC-xi2}
\E[\xi_{n+1}^2|\cF_n]=1 -a^2
\end{equation}
and finally
\begin{align}
\label{ESPC-xi-eps}
\E[\eps_{n+1}\xi_{n+1}|\cF_n] &= \E\big[\big((1-\gamma_n)Y_n+(\alpha_{n+1}+c)X_{\beta_{n+1}}\big)(\alpha_{n+1}-a)X_{\beta_{n+1}}|\cF_n\big]\nonumber \\
		&= \E\big[\big((1-\gamma_n)(\alpha_{n+1}-a)Y_nX_{\beta_{n+1}}+(\alpha_{n+1}+c)(\alpha_{n+1}-a)|\cF_n\big]\nonumber \\
		&= 1-a^2.
\end{align}
Hereafter, we deduce from \eqref{quad-Mn}, \eqref{ESPC-eps2}, \eqref{ESPC-xi2} and \eqref{ESPC-xi-eps} that
\begin{align*}
\E\big[(\Delta\cM_{n+1})(\Delta\cM_{n+1})^T|\cF_n\big]&= 
a_{n+1}^2 \big(1 +2ac+c^2 - (\gamma_k-1)^2Y_k^2\big)\begin{pmatrix}0 & 0\\ 0 & 1 \end{pmatrix} \\ 
& \quad + a_{n+1}\frac{c}{a+c}(1 -a^2)\begin{pmatrix}0 & 1\\ 1 & 0 \end{pmatrix} \\ 
& \quad + \Big(\frac{c}{a+c}\Big)^2(1 -a^2)\begin{pmatrix}1 & 0\\ 0 & 0 \end{pmatrix}.
\end{align*}
We are now able to compute the quadratic variation of $\cM_n$, that is
\begin{align*}
\langle \cM\rangle_n & = \sum_{k=0}^{n-1} 
a_{k+1}^2 \big(1 +2ac+c^2 - (\gamma_k-1)^2Y_k^2\big)\begin{pmatrix}0 & 0\\ 0 & 1 \end{pmatrix} \nonumber\\ 
& \quad + \sum_{k=0}^{n-1} a_{k+1}\frac{c}{a+c}(1 -a^2)\begin{pmatrix}0 & 1\\ 1 & 0 \end{pmatrix} \nonumber\\ 
& \quad + n\Big(\frac{c}{a+c}\Big)^2(1 -a^2)\begin{pmatrix}1 & 0\\ 0 & 0 \end{pmatrix}.
\end{align*}
Consequently, 
\begin{align}
\label{QV-MMN}
\langle \cM\rangle_n &= v_n(1 +2ac+c^2)\begin{pmatrix}0 & 0\\ 0 & 1 \end{pmatrix} + \nonumber
w_n\frac{c}{a+c}(1 -a^2)\begin{pmatrix}0 & 1\\ 1 & 0 \end{pmatrix} 
\\ & \quad + n\Big(\frac{c}{a+c}\Big)^2(1 -a^2)\begin{pmatrix}1 & 0\\ 0 & 0 \end{pmatrix} -\cR_n\begin{pmatrix}0 & 0\\ 0 & 1 \end{pmatrix}
\end{align}
where
\begin{equation*}
\label{DEF-sumVN}
v_n = \sum_{k=1}^n a_k^2 \ , \quad w_n=\sum_{k=1}^n a_k \quad\text{and}\quad
	\cR_n= \sum_{k=0}^{n-1}a_{k+1}^2(\gamma_k-1)^2Y_k^2.
\end{equation*}
Hereafter, we immediately deduce from \eqref{QV-MMN} that 
\begin{equation}
\label{QV-M}
\cMc_n = (1 +2ac+c^2)\sum_{k=1}^n a_k^2 - \cR_n 
\end{equation}
and that
\begin{equation}
\label{QV-N}
\langle N\rangle_n = \Big(\frac{c}{a+c}\Big)^2(1 -a^2)n.
\end{equation}
The asympotic behavior of $M_n$ is closely related to the one of $(v_n)$
as one can observe that we always have $\cMc_n \leq (1 +2ac+c^2) v_n$ and thus $\cMc_n = O(v_n)$. Consequently to the definition of $(a_n)$, we have three regimes of behavior for $(M_n)$. In the diffusive regime where $a<(1-c)/2$,
\begin{equation}
\label{VN-DIF}
\lim_{n\to\infty} \frac{v_n}{n^{1-2(a+c)\lambda}}=\ell \quad\text{where}\quad
\ell=\frac{1}{1-2(a+c)\lambda}\Big(\frac{\Gamma(1+a\lambda)}{\Gamma(\lambda)}\Big)^2.
\end{equation}
In the critical regime where $a=(1-c)/2$,
\begin{equation}
\label{VN-CR}
\lim_{n\to\infty} \frac{v_n}{\log n}=\Big(\frac{\Gamma(\frac{c+3}{2(c+1)})}{\Gamma(\frac{1}{c+1})}\Big)^2.
\end{equation}
In the superdiffusive regime where $a>(1-c)/2$,
\begin{equation}
\label{VN-SDIF}
\lim_{n\to\infty} {v_n}= \sum_{n=1}^\infty \Big(\frac{\Gamma(n-c\lambda)\Gamma(1+a\lambda)}{\Gamma(n+a\lambda)\Gamma(\lambda)}\Big)^2.
\end{equation}

%
\section{Proofs of the almost sure convergence results}
\label{S-AS-results}

\begin{lemma}
\label{lemMV} 
Let $(V_n)$ be the sequence of positive definite diagonal matrices of order $2$ given by
\begin{equation}
\label{DEF-VN}
 V_n = \frac{1}{\sqrt{n}}\begin{pmatrix} 1 & 0 \\ 0 &\dfrac{a}{a+c}a_n^{-1}\end{pmatrix}.
\end{equation}
Then, the quadratric variation of $\langle \cM\rangle_n$ satisfies in the diffusive regime where $a<(1-c)/2$,
\begin{equation}
\label{lim-VMM}
\lim_{n\to\infty} V_n \langle \cM\rangle_n V_n = V \quad\text{a.s.}
\end{equation}
where the matrix $V$ is given by
\begin{equation}
\label{DEF-V}
V = \frac{1}{(a+c)^2}\begin{pmatrix}c^2(1 -a^2) & ac(c+1)(1 +a)
\\[6pt] ac(c+1)(1 +a) & \dfrac{a^2(1 +2ac+c^2)(c+1)}{1-c-2a} \end{pmatrix}.
\end{equation}
\end{lemma}

\begin{remark}
\label{remWn}
Following the same steps as in the proof of Lemma \ref{lemMV}, we find that in the critical regime $a = (1-c)/2$, the sequence of normalization matrices $(V_n)$ has to be replaced by 
\begin{equation}
\label{DEF-WN}
W_n = \frac{1}{\sqrt{n\log n}}\begin{pmatrix} 1 & 0 \\[5pt] 0 & \dfrac{a}{a+c}a_n^{-1}\end{pmatrix}.
\end{equation}
The limit matrix $V$ also need to be replaced by 
\begin{equation}
\label{DEF-W}
W = \frac{(c-1)^2}{c+1}\begin{pmatrix} 0 & 0 \\ 0 & 1\end{pmatrix}.
\end{equation}
\end{remark}
\begin{proof}{Lemma}{\ref{lemMV}}
We immediately obtain from Theorem \ref{THM-CVPS-DIF} and \eqref{lim-CN}, \eqref{QV-MMN}, \eqref{VN-DIF} that
\begin{align*}
\lim_{n\to\infty} V_n \langle \cM\rangle_n V_n^T
& =\Big(\frac{a}{a+c}\Big)^2\frac{1}{1-2\lambda(a+c)}(1 +2ac+c^2)\begin{pmatrix}0 & 0\\ 0 & 1 \end{pmatrix} \\
& \quad + \frac{1}{1-\lambda}(1 -a^2)\frac{ac}{(a+c)^2}\begin{pmatrix}0 & 1 \\ 1 & 0 \end{pmatrix} \\
& \quad + (1 -a^2)\Big(\frac{c}{a+c}\Big)^2\begin{pmatrix}1 & 0\\ 0 & 0 \end{pmatrix}\\
& = \frac{1}{(a+c)^2}\begin{pmatrix}c^2(1 -a^2) & ac(c+1)(1 +a)
\\[6pt] ac(c+1)(1 +a) & \dfrac{a^2(1 +2ac+c^2)(c+1)}{1-c-2a} \end{pmatrix}
\end{align*}
which is exactly what we wanted to prove.
\end{proof}

\subsection{The diffusive regime}
%
\begin{proof}{Theorem}{\ref{THM-CVPS-DIF}}
We shall make extensive use of the strong law of large numbers for martingales given, e.g. by theorem 1.3.24 of \cite{Duflo1997}. First, we have for $M_n$ that for any $\gamma>0$,
\begin{equation*}
M_n^2 = O\big((\log v_n)^{1+\gamma}v_n\big)\quad\text{a.s.}
\end{equation*}
Then, by definition of $M_n$ and as $a_n$ is asymptotically equivalent to $n^{-(a+c)\lambda}$ and $v_n$ is asymptotically equivalent to $n^{1-2(a+c)\lambda}$, it ensures
\begin{equation*}
\frac{Y_n^2}{n^2} = O\Big((\log n)^{1+\gamma}\frac{n^{1-2(a+c)\lambda}}{n^{2(1-(a+c)\lambda)}}\Big)\quad\text{a.s.}
\end{equation*}
and finally 
\begin{equation*}
\frac{Y_n^2}{n^2} = O\Big(\frac{(\log n)^{1+\gamma}}{n}\Big)\quad\text{a.s.}
\end{equation*}
This implies that 
\begin{equation}
\label{lim-YN}
\lim_{n\to\infty} \frac{Y_n}{n}=0 \quad\text{a.s.}
\end{equation}
We now focus our attention on $N_n$. By the same token as before, we have that for any $\gamma>0$,
\begin{equation*}
N_n^2 = O\big((\log n)^{1+\gamma}n\big)\quad\text{a.s.}
\end{equation*}
which by definition of $N_n$ gives us
\begin{equation*}
\frac{\big(S_n-\frac{a}{a+c}Y_n\big)^2}{n^2} = O\Big(\frac{(\log n)^{1+\gamma}}{n}\Big)\quad\text{a.s.}
\end{equation*}
and we conclude
\begin{equation}
\label{lim-SYN}
\lim_{n\to\infty} \frac{S_n}{n}-\frac{a}{a+c}\frac{Y_n}{n}=0 \quad\text{a.s.}
\end{equation}
This achieves the proof of Theorem \ref{THM-CVPS-DIF} as the convergences \eqref{lim-YN} and \eqref{lim-SYN} hold almost surely.
\end{proof}

\begin{proof}{Theorem}{\ref{THM-LFQ-DIF}}
We need to check that all the hypotheses of Theorem A.2 in \cite{LaulinBercu2020} are satisfied. Thanks to Lemma \ref{lemMV}, hypothesis $(\textnormal{H.1})$ holds almost surely. In order to verify that Lindeberg's condition $(\textnormal{H.2})$ is satisfied, we have from \eqref{DEF-M-N} together with \eqref{DEF-MMN} and $V_n$ given by \eqref{DEF-VN}
that for all $1 \leq k \leq n$
\begin{equation*}
V_n \Delta \cM_k =\frac{1}{(a+c)\sqrt{n}}
\begin{pmatrix} c\xi_{n+1}\\ a a_n^{-1}a_{k}\eps_{k}\end{pmatrix}
\end{equation*}
which implies that 
\begin{equation}
\label{VNMMN-BOUND2}
\| V_n \Delta \cM_k \|^2 
= \frac{1}{(a+c)^2n}  \big(c^2+ a^2 a_n^{-2}a_k^2 \eps_k^2\big)
\end{equation}
and
\begin{equation}
\label{VNMMN-BOUND4}
\| V_n \Delta \cM_k \|^4 
= \frac{1}{(a+c)^4n^2}  \big(c^4 +2a^2c^2a_n^{-2}a_k^2 \eps_k^2+a_n^{-4}a_k^4\eps_k^4\big).
\end{equation}
Consequently, we obtain that for all $\eps >0$,
\begin{equation}
\sum_{k=1}^n \E\big[\|V_n \Delta \cM_k \|^2 \1_{\{\|V_n\Delta \cM_k \|>\eps\}}\big|\cF_{k-1}\big]
    \leq  \frac{1}{\eps^2}\sum_{k=1}^n \E\big[\|V_n \Delta \cM_k \|^4\big|\cF_{k-1}\big]. 
\label{LINDEBERG-DR1} 
\end{equation}
It follows from \eqref{lim-CN} that
\begin{align*}
a_n^{-2}\sum_{k=1}^n a_k^2  = O(n) \quad\text{and}\quad
a_n^{-4}\sum_{k=1}^n a_k^4  = O(n).
\end{align*}
Hence, using that the sequence $(\eps_n)$ is uniformly bounded
\begin{equation}
\label{SUPEPS}
\sup_{1\leq k\leq n} |\eps_k| \leq c+2 \quad\text{a.s.}
\end{equation}
we find that
\begin{equation*}
\sum_{k=1}^n \E\big[\|V_n \Delta \cM_k \|^4 \big|\cF_{k-1}\big] = O\Big(\frac{1}{n}\Big)\quad\text{a.s.}
\end{equation*}
which ensures that Lindeberg's condition $(\textnormal{H.2})$ holds almost surely, that is
for all $\eps >0$,
\begin{equation}
\label{LINDEBERG-DR4}
\lim_{n \rightarrow \infty} \sum_{k=1}^n \E\big[\|V_n \Delta \cM_k \|^2 \1_{\{\|V_n\Delta \cM_k \|>\eps\}}\big|\cF_{k-1}\big]= 0 \quad \text{a.s.}
\end{equation}
Hereafter, we need to verify $(\textnormal{H.3})$ is satisfied in the special case $\beta=2$ that is
\begin{equation*}
\label{LFQ-DIF1}
\sum_{n=1}^{\infty} \frac{1}{\bigl(\log  (\det V_{n}^{-1})^2\bigr)^{2}}\E\big[\|V_{n} \Delta \cM_{n}\|^{4}\big|\cF_{n-1}\big]<\infty \quad \text{a.s.}
\end{equation*}
We immediately have from \eqref{DEF-VN}
\begin{equation}
\label{DET-VN}
\det V_n^{-1} = \frac{a+c}{a}\sqrt{n}a_n.
\end{equation} 
Hence, we obtain from \eqref{lim-CN} and \eqref{DET-VN} that
\begin{equation}
\label{LOG-DET-VN}
\lim_{n \rightarrow \infty}  \frac{\log  (\det V_{n}^{-1})^2 }{\log n} = 1-2(a+c)\lambda.
\end{equation}
Therefore, we can replace $\log  (\det V_{n}^{-1})^2$ by $\log n$ in \eqref{LFQ-DIF1}. Hereafter, we obtain from \eqref{VNMMN-BOUND4} and \eqref{SUPEPS} that
\begin{align}
\label{LFQ-DIF2}
\sum_{n=2}^{\infty} \frac{1}{(\log n)^2}\E\big[\|V_{n} \Delta \cM_{n}\|^{4}\big|\cF_{n-1}\big]
& = O \Big( \sum_{n=1}^{\infty} \frac{1}{(n\log n)^2}  \Big). 
\end{align}
Thus, \eqref{LFQ-DIF2} guarentees that $(\textnormal{H.3})$ is verified.
We are now going to apply the quadratic strong law given by Theorem A.2 in \cite{LaulinBercu2020}. We get from equation \eqref{LOG-DET-VN} that
\begin{equation}
\label{LFQ-DIF3}
\lim_{n \rightarrow \infty}  \frac{1}{\log n}\sum_{k=1}^n
\Big(\frac{(\det V_{k})^2 - (\det V_{k+1})^2}{(\det V_k)^2}\Big)V_k\cM_k\cM_k^T V_k^T =
(1-2(a+c)\lambda) V \quad \text{a.s.}
\end{equation}
However, we obtain from \eqref{lim-CN} and \eqref{DET-VN} that
\begin{equation}
\label{lim-DET-VN}
\lim_{n \rightarrow \infty} n\Big(\frac{(\det V_{n})^2 - (\det V_{n+1})^2}{(\det V_n)^2}\Big)=1-2(a+c)\lambda.
\end{equation}
Finally, let $u = \begin{pmatrix}1, 1\end{pmatrix}^T$ we have
\begin{equation}
\label{SN-vVNMMN}
u^T V_n\cM_n = \frac{S_n}{\sqrt{n}}
\end{equation} 
and we deduce from  \eqref{LFQ-DIF3}, \eqref{lim-DET-VN} and \eqref{SN-vVNMMN} that
\begin{equation}
\label{LFQ-DIF4}
\lim_{n \rightarrow \infty}  \frac{1}{\log n}\sum_{k=1}^n\frac{S_k^2}{k^2} = (1-2(a+c)\lambda) v^TVv \quad \text{a.s.}
\end{equation}
which, together with
\begin{equation}
\label{COV-vVv}
u^TVu = \frac{2ac+c-1}{2a+c-1}
\end{equation}
completes the proof of Theorem \ref{THM-LFQ-DIF}.
\end{proof}

\subsection{The critical regime}
%
\begin{proof}{Theorem}{\ref{THM-CVPS-CR}}
Again, we shall make use of the strong law of large numbers for martingales given, e.g. by theorem 1.3.24 of \cite{Duflo1997}. First, we have for $M_n$ that for any $\gamma>0$,
\begin{equation*}
M_n^2 = O\big((\log v_n)^{1+\gamma}v_n\big)\quad\text{a.s.}
\end{equation*}
which by definition of $M_n$ and as $a_n$ is asymptotically equivalent to $n^{-1/2}$ and $v_n$ is asymptotically equivalent to $\log n$ ensures that
\begin{equation*}
\frac{Y_n^2}{(\sqrt{n}\log n)^2} = O\Big((\log \log n)^{1+\gamma} \frac{\log n}{(\log n)^2} \Big)\quad\text{a.s.}
\end{equation*}
and finally that
\begin{equation*}
\frac{Y_n^2}{(\sqrt{n}\log n)^2} = O\Big(\frac{(\log \log n)^{1+\gamma}}{\log n}\Big)\quad\text{a.s.}
\end{equation*}
This implies that 
\begin{equation}
\label{lim-YN-CR}
\lim_{n\to\infty} \frac{Y_n}{\sqrt{n}\log n}=0 \quad\text{a.s.}
\end{equation}
In addition, we still have that for any $\gamma>0$,
\begin{equation*}
N_n^2 = O\big((\log n)^{1+\gamma}n\big)\quad\text{a.s.}
\end{equation*}
which by definition of $N_n$ gives us
\begin{equation*}
\frac{\big(S_n-\frac{a}{a+c}Y_n\big)^2}{(\sqrt{n}\log n)^2} = O\Big((\log n)^{\gamma-1}\Big)\quad\text{a.s.}
\end{equation*}
Taking e.g. $\gamma=\frac{1}{2}$ we can conclude that
\begin{equation}
\label{lim-SYN-CR}
\lim_{n\to\infty} \frac{S_n}{\sqrt{n}\log n}-\frac{a}{a+c}\frac{Y_n}{\sqrt{n}\log n}=0 \quad\text{a.s.}
\end{equation}
This achieves the proof of Theorem \ref{THM-CVPS-CR} as the convergences \eqref{lim-YN-CR} and \eqref{lim-SYN-CR} hold almost surely.
\end{proof}
{}
\begin{proof}{Theorem}{\ref{THM-LFQ-LIL-CR}}
The proof of the quadratic strong law \eqref{QUAD-CR} is left to the reader as it 
follows essentially the same lines as that of \eqref{QUAD-DIF}. The only minor change is that the matrix $V_n$ has to be replaced by the matrix $W_n$ defined in \eqref{DEF-WN}.
We shall now proceed to the proof of the law of iterated logarithm
given by \eqref{LIL-CR}. On the one hand, it follows from \eqref{lim-CN} and \eqref{VN-DIF} that 
\begin{equation}
\label{CONDLIL-CR}
\sum_{n=1}^{+\infty} \frac{a_n^4}{v_n^2} < \infty.
\end{equation}
Moreover, we have from \eqref{QV-M} and \eqref{QV-N} that 
\begin{equation*}
\lim_{n\to\infty}\frac{\langle M\rangle_n}{v_n} =  1 +2ac+c^2\quad\text{a.s.}
\quad\text{and}\quad
\lim_{n\to\infty}\frac{\langle N\rangle_n}{n} = \Big(\frac{c}{a+c}\Big)^2(1 -a^2) \quad\text{a.s.}
\end{equation*}
Consequently, we deduce from the law of iterated logarithm for martingales due to Stout \cite{Stout1974}, see also Corollary 6.4.25 in \cite{Duflo1997},  
that $(M_n)$ satisfies when $a=(1-c)/2$
\begin{align*}
\limsup_{n \rightarrow \infty} \frac{M_n}{(2 v_n \log \log v_n)^{1/2}}  & = -\liminf_{n \rightarrow \infty} \frac{M_n}{(2 v_n \log \log v_n)^{1/2}}  \\
& = \sqrt{1 +c}  \quad\text{a.s.}
\end{align*}
However, as $a_n v_n^{-1/2}$ is asymptotically equivalent to $(n\log n)^{-1/2}$, we immediately obtain from \eqref{VN-CR} that
\begin{align}
\label{LIL-CR1}
\limsup_{n \rightarrow \infty} \frac{Y_n}{(2 n \log n \log \log \log n)^{1/2}}  & = -\liminf_{n \rightarrow \infty} \frac{Y_n}{(2 n \log n \log \log \log n)^{1/2}} \nonumber  \\
& = \sqrt{1 +c}  \quad\text{a.s.}
\end{align}
The law of iterated logarithm for martingales also allow us to find that $(N_n)$ satisfies
\begin{align*}
\label{LIL-CR2}
\limsup_{n \rightarrow \infty} \frac{N_n}{(2 n \log \log n)^{1/2}}  & = -\liminf_{n \rightarrow \infty} \frac{N_n}{(2 n \log \log n)^{1/2}}  \\
& = \frac{2c}{c+1}\sqrt{(1 -a^2)}  \quad\text{a.s.}
\end{align*}
which ensures that
\begin{equation*}
\limsup_{n \rightarrow \infty}  \frac{N_n}{(2 n \log n \log \log \log n)^{1/2}} 
= 0 \quad\text{a.s.} 
\end{equation*} 
Hence, we deduce from \eqref{SN-MG} and \eqref{LIL-CR1} that
\begin{align*}
\label{LIL-CR3}
\limsup_{n \rightarrow \infty} \frac{S_n}{(2 n\log n \log \log \log n)^{1/2}}  
& = \limsup_{n \rightarrow \infty} \frac{N_n +\frac{1-c}{1+c}a_n^{-1}M_n}{(2 n\log n \log \log \log n)^{1/2}} \nonumber \\
& = \limsup_{n \rightarrow \infty} \frac{1-c}{1+c}\frac{Y_n}{(2 n\log n \log \log \log n)^{1/2}} \nonumber \\
& = -\liminf_{n \rightarrow \infty} \frac{1-c}{1+c}\frac{Y_n}{(2 n\log n \log \log \log n)^{1/2}} \nonumber  \\
& = -\liminf_{n \rightarrow \infty} \frac{S_n}{(2 n\log n \log \log \log n)^{1/2}}. \nonumber
\end{align*}
Hence, we obtain that 
\begin{align*}
\limsup_{n \rightarrow \infty} \frac{S_n^2}{2 n\log n \log \log \log n}
&= \limsup_{n \rightarrow \infty} \Big(\frac{1-c}{1+c}\Big)^2\frac{Y_n^2}{2 n\log n \log \log \log n} \\
&=\frac{(1-c)^2}{1+c}
\end{align*}
which immediately leads to \eqref{LIL-CR}, thus completing the proof of Theorem \ref{THM-LFQ-LIL-CR}.
\end{proof}

\subsection{Superdiffusive regime}
%
\begin{proof}{Theorem}{\ref{THM-CVPS-SDIF}}
Hereafter, we shall again make extensive use of the strong law of large numbers for martingales given, e.g. by theorem 1.3.24 of \cite{Duflo1997} in order to prove \eqref{CVPS-SDIF}. When $a>(1-c)/2$, we have from \eqref{VN-SDIF} that $v_n$ converges. Hence, as $\langle M\rangle_n \leq (1+2ac+c^2)v_n$, we clealy have that $\langle M\rangle_\infty <\infty$ almost surely and we can conclude that
\begin{equation*}
\lim_{n\to\infty} M_n = M \quad\text{a.s.} \quad\text{where}\quad M = \sum_{k=1}^{\infty} a_k \eps_k
\end{equation*}
which by definition of $M_n$ and as $a_n$ is asymptotically equivalent to $\frac{\Gamma(1+a\lambda)}{\Gamma(\lambda)}n^{-(a+c)\lambda}$ ensures that
\begin{equation}
\label{lim-YN-SDIF}
\lim_{n\to\infty}\frac{Y_n}{n^{(a+c)\lambda}} = Y \quad\text{a.s.}\quad\text{where}\quad Y = \frac{\Gamma(\lambda)}{\Gamma(1+a\lambda)}M.
\end{equation}
Moreover, we still have that for any $\gamma>0$,
\begin{equation*}
N_n^2 = O\big((\log n)^{1+\gamma}n\big)\quad\text{a.s.}
\end{equation*}
which by definition of $N_n$ gives us for all $t\geq 0$
\begin{equation*}
\frac{\big(S_n-\frac{a}{a+c}Y_n\big)^2}{n^{2(a+c)\lambda}} = O\Big(\frac{(\log n)^{1+\gamma}}{n^{2(a+c)\lambda-1}}\Big)\quad\text{a.s.}
\end{equation*}
As $a>(1-c)/2$ in the superdiffusive regime, we obtain thanks to \eqref{lim-YN} that for all $t\geq 0$
\begin{equation}
\label{lim-SYN-SDIF1}
\lim_{n\to\infty} \frac{S_{\floor{nt}}}{{\floor{nt}}^{(a+c)\lambda}}-\frac{a}{a+c}\frac{Y_{\floor{nt}}}{{\floor{nt}}^{(a+c)\lambda}}=0 \quad\text{a.s.}
\end{equation}
The convergences \eqref{lim-YN-SDIF} and \eqref{lim-SYN-SDIF1} hold almost surely and $\floor{nt}$ is asymptotically equivalent to $nt$ which implies
\begin{equation}
\label{lim-SYN-SDIF-2}
\lim_{n\to\infty} \frac{S_{\floor{nt}}}{n^{(a+c)\lambda}}= t^{(a+c)\lambda}L_c \quad\text{a.s.}
\end{equation} 
Finally, the fact that \eqref{lim-SYN-SDIF-2} holds almost surely ensures that it also holds for the finite-dimensional distributions, and we obtain \eqref{CVPS-cont-SDIF} with $\Lambda_t = t^{(a+c)\lambda}L_c$ and $L_c =\frac{a}{a+c}Y$. \medskip

We shall now proceed to the proof of the mean square convergence \eqref{CVL2-SDIF}. On the one hand, as $M_0=0$ we have from \eqref{QV-M} that
\begin{equation*}
\E\big[M_n^2\big]=\E\big[\langle M\rangle_n\big]\leq (1+2ac+c^2)v_n.
\end{equation*}
Hence, we obtain from \eqref{VN-SDIF} that
\begin{equation*}
\sup_{n\geq 1} \E\big[M_n^2\big] < \infty
\end{equation*}
which ensures that the martingale $(M_n)$ is bounded in $\IL^2$. Therefore, we have the mean square convergence
\begin{equation*}
\lim_{n\to\infty} \E\big[\big|M_n-M\big|^2\big]=0
\end{equation*}
which implies that
\begin{equation}
\label{L2-Y}
\lim_{n\to\infty} \E\Big[\Big|\frac{Y_n}{n^{(a+c)\lambda}}-Y\Big|^2\Big]=0.
\end{equation}
On the other hand, for any $n\geq0$, the martingale $(N_n)$ satisfies
\begin{equation*}
\E\big[N_n^2\big]=\E\big[\langle N\rangle_n\big]\leq (1 -a^2)\Big(\frac{c}{a+c}\Big)^2n
\end{equation*}
and since $(a+c)\lambda>\frac{1}{2}$ we obtain
\begin{equation}
\label{L2-N}
\lim_{n\to\infty} \E\Big[\Big|\frac{N_n}{n^{(a+c)\lambda}}\Big|^2\Big]=0.
\end{equation}
Finally, we obtain the mean square convergence \eqref{CVL2-SDIF} from \eqref{L2-Y} and \eqref{L2-N} and we achieve the proof of Theorem \ref{THM-CVPS-SDIF}.
\end{proof}

\begin{proof}{Theorem}{\ref{THM-L2-SDIF}}
We start by the calculation of the expectation \eqref{ESP-SDIF}. We immediately have from \eqref{ESPC-Y} that
\begin{equation*}
\E[Y_{n+1}]=\gamma_n \E[Y_n] = \Big(\frac{n+a\lambda}{n-c\lambda}\Big)\E[Y_n]
\end{equation*}
which leads to
\begin{equation}
\label{EQ-ESP-Y}
\E[Y_n]=\prod_{k=1}^{n-1}\Big(\frac{k+a\lambda}{k-c\lambda}\Big)\E[Y_1] = \prod_{k=1}^{n-1}\Big(\frac{k+a\lambda}{k-c\lambda}\Big)\E[X_1]
=(2q-1)a_n^{-1}.
\end{equation}
Hence, we immediately get equation \eqref{ESP-SDIF} from \eqref{EQ-ESP-Y}, that is
\begin{equation*}
\E[L_c] = \frac{a\Gamma(\lambda)}{(a+c)\Gamma(1+a\lambda)}\E[M]=\frac{a\Gamma(\lambda)}{(a+c)\Gamma(1+a\lambda)}\E[M_n]= \frac{a(2q-1)\Gamma(\lambda)}{(a+c)\Gamma(1+a\lambda)}.
\end{equation*}
Hereafter, we obtain from \eqref{ESPC-eps2} by taking expectation on both sides that
\begin{equation*}
\E[Y_{n+1}^2]=1+2ac+c^2 +(2\gamma_n-1)\E[Y_n^2]= 
1 +2ac +c^2 +\Big(\frac{n+(2a+c)\lambda}{n-c\lambda}\Big)\E[Y_n^2]
\end{equation*}
and thanks to well-kown recursive relation solutions and Lemma B.1 in \cite{Bercu2018}, we get
\begin{align*}
\E[Y_n^2]& = (1+2ac+c^2)\prod_{k=0}^{n-1}\Big(\frac{k+(2a+c)\lambda}{k-c\lambda}\Big)
		\sum_{k=0}^{n-1}\prod_{i=0}^{k}\frac{i-c\lambda}{i+(2a+c)\lambda} \\
		& = \frac{(1+2ac+c^2)\Gamma(n+(2a+c)\lambda)\Gamma(\lambda)}{\Gamma(n-c\lambda)\Gamma(1+(2a+c)\lambda)}\sum_{k=0}^{n-1}\frac{\Gamma(k+\lambda)\Gamma(1+(2a+c)\lambda)}{\Gamma(k+1+(2a+c)\lambda)\Gamma(\lambda)}\\
		& = \frac{(1+2ac+c^2)\Gamma(n+(2a+c)\lambda)}{\Gamma(n-c\lambda)}\sum_{k=1}^{n}\frac{\Gamma(k+\lambda-1)}{\Gamma(k+(2a+c)\lambda)}\\
		& = \frac{(1+2ac+c^2)\Gamma(n+(2a+c)\lambda)}{\lambda(2a+c-1)\Gamma(n-c\lambda)}
		\Big(\frac{\Gamma(\lambda)}{\Gamma((2a+c)\lambda)}
		-\frac{\Gamma(n+\lambda)}{\Gamma(n+(2a+c)\lambda)}\Big).
\end{align*}
Hence, we obtain from \eqref{lim-CN}, \eqref{DEF-M-N} and \eqref{L2-Y} that
\begin{equation}
\E[Y^2] = \lim_{n\to\infty} \frac{\E[Y_n^2]}{n^{2(a+c)\lambda}} =\frac{(1+2ac+c^2)\Gamma(\lambda)}{\lambda(2a+c-1)\Gamma((2a+c)\lambda)}.
\end{equation}
\end{proof}


\section{\blue{Proofs of the functional limit theorems}}
\label{S-FCLT-results}

\subsection{The diffusive regime}
%

\begin{proof}{Theorem}{\ref{CLFT-RERW-DR}}
In order to apply Theorem \ref{T-CLTF} in the Appendix, we must verify that (H.1), (H.2) and (H.3) are satisfied. \medskip \\
\textbf{(H.1)} We have from \eqref{lim-VMM} and the fact that $a_{\floor{nt}}$ is asymtotically equivalent to $t^{-(a+c)\lambda}a_n$ that
\begin{equation*}
V_n \langle \cM \rangle_{\floor{nt}} V_n^T {\underset{n\to\infty}{\longrightarrow}} V_t \quad\text{a.s.}
\end{equation*}
where
\begin{equation*}
V_t = \frac{1}{(a+c)^2}\begin{pmatrix}c^2(1 -a^2)t & ac(c+1)(1 +a)t^{1-(a+c)\lambda}
\\[6pt] ac(c+1)(1 +a)t^{1-(a+c)\lambda} & \dfrac{a^2(1 +2ac+c^2)(c+1)}{1-c-2a}t^{1-2(a+c)\lambda} \end{pmatrix}.
\end{equation*}
\medskip \\
\textbf{(H.2)} We also get that Lindeberg's condition is satisfied as we already know from \eqref{LINDEBERG-DR4} that for all $\eps>0$
\begin{equation*}
\lim_{n \rightarrow \infty} \sum_{k=1}^n \E\big[\|V_n \Delta \cM_k \|^2 \1_{\{\|V_n\Delta \cM_k \|>\eps\}}\big|\cF_{k-1}\big]= 0 \quad \text{a.s.}
\end{equation*}
which implies from \eqref{VNMMN-BOUND4} and the fact that $V_nV_{\floor{nt}}^{-1}$ converges
\begin{align*}
\lim_{n\to\infty}\sum_{k=1}^{\floor{nt}} \E\bigl[\|V_n \Delta \cM_k \|^2 \1_{\{\|V_n\Delta \cM_k \|>\eps\}}\bigl|\cF_{k-1}\bigr] & \leq \lim_{n\to\infty}\sum_{k=1}^{\floor{nt}} \E\bigl[\|V_{n} \Delta \cM_k \|^4 \bigr] \\
& \leq\lim_{n\to\infty}\sum_{k=1}^{\floor{nt}} \E\bigl[\|(V_nV_{\floor{nt}}^{-1})V_{\floor{nt}} \Delta \cM_k \|^4\bigr] \\
& = 0\quad\text{a.s.}
\end{align*}
\medskip \\
\textbf{(H.3)} In this particular case, we have $V_t = t K_1 + t^{\alpha_2} K_2 + t^{\alpha_3} K_3$ where 
\begin{equation*}
{\alpha_2} = {1-(a+c)\lambda}>0 \quad\text{and}\quad {\alpha_3}={1-2(a+c)\lambda}>0
\end{equation*} 
as $a\leq(1-c)/2$, and the matrix are symmetric
\begin{equation*}
K_1 =\frac{c^2(1-a^2)}{(a+c)^2}\begin{pmatrix}1 & 0 \\ 0 & 0\end{pmatrix}, \
K_2 =\frac{ac(c+1)(a+1)}{(a+c)^2}\begin{pmatrix}0 & 1 \\ 1 & 0\end{pmatrix},
\end{equation*}
\begin{equation*}
K_3 =\frac{a^2(1+2ac+c^2)(c+1)}{(1-2a-c)(a+c)^2}\begin{pmatrix}0 & 0 \\ 0 & 1\end{pmatrix}.
\end{equation*}
\medskip \\
Consequently, we obtain that
\begin{equation*}
\label{CLTF-MERW}
\big(V_n \cM_{\floor{nt}}, \ {t \geq 0}\big) \Longrightarrow \big(\cB_{t}, \ {t \geq 0}\big)
\end{equation*}
where $\cB$ is defined as in \eqref{CLTF}.
Finally,  using the fact that $S_{\floor{nt}}$ is asymptotically equivalent to $N_{\floor{nt}}+t^{(a+c)\lambda}\frac{a}{a+c}a_{n}^{-1}M_{\floor{nt}}$ and  
 multiplying $u_t=\comb{1}{t^{(a+c)\lambda}}$ we conclude
\begin{equation}
\label{CLTF-MERW}
\big(\frac{1}{\sqrt{n}} S_{\floor{nt}}, \ {t \geq 0}\big) \Longrightarrow \big(W_t, \ {t \geq 0}\big)
\end{equation}
where $W_t =u_t^T \cB_t$. It only remains to compute the covariance function of $W$ that is for $0\leq s \leq t$
\begin{align*}
\E\big[W_sW_t\big] & = u_s^T\E\big[\cB_s \cB_t^T\big]u_t \\
& = u_s^TV_su_t \\
& = u_s^T\big(sK_1+ s^{1-(a+c)\lambda}K_2+s^{1-2(a+c)\lambda}K_3)u_t \\
& = \frac{c^2(1-a^2)}{(a+c)^2}s +\frac{ac(c+1)(a+1)}{(a+c)^2}s^{1-(a+c)\lambda}(s^{(a+c)\lambda}+t^{(a+c)\lambda})\\
&\quad +\frac{a^2(1+2ac+c^2)(c+1)}{(1-2a-c)(a+c)^2}s^{1-2(a+c)\lambda}(st)^{(a+c)\lambda} \\
& = \Big(\frac{c^2(1-a^2)}{(a+c)^2} +\frac{ac(c+1)(a+1)}{(a+c)^2}\Big) s \\
&\quad + \Big(\frac{ac(c+1)(a+1)}{(a+c)^2}+\frac{a^2(1+2ac+c^2)(c+1)}{(1-2a-c)(a+c)^2}\Big)s\Big(\frac{t}{s}\Big)^{(a+c)\lambda}\\
& = \frac{c(a+1)}{a+c}s + \frac{a(1-c^2)}{(a+c)(1-2a-c)}s\Big(\frac{t}{s}\Big)^{(a+c)\lambda}.
\end{align*}
\end{proof}

\begin{proof}{Corollary}{\ref{coro-Gn}}
As for Corollary 4.1 from \cite{Bertenghi2020}, we observe that
\begin{equation*}
\frac{G_n}{\sqrt{n}} = \int_0^1 \frac{S_{\floor{nt}}}{\sqrt{n}}\drm t.
\end{equation*}
Consquently, $G_n/\sqrt{n}$ is a continuous function of ${S_{\floor{nt}}}/{\sqrt{n}}$ in $D([0,1])$. Hence, the functional distribution from Theorem \ref{CLFT-RERW-DR} gives us that
\begin{equation*}
\frac{G_n}{\sqrt{n}} = \int_0^1 \frac{S_{\floor{nt}}}{\sqrt{n}}\drm t \liml \int_0^1 W_t \drm t.
\end{equation*}
The process $\big(W_t, \ t\geq0 \big)$ is a continuous real-valued and centered Gaussian process starting from the origin, which implies that $\int_0^1 W_t \drm t$ is also one. Its covariance is given by
\begin{align*}
\E\Big[\big(\int_0^1 W_s \drm s\big)\big(\int_0^1 W_t \drm t\big)\Big] & = 2\int_0^1\int_0^t	\E\big[ W_s W_t\big] \drm s \drm t \\
& = 2\frac{a(1-c^2)}{(a+c)(1-2a-c)}\int_0^1\int_0^t s\Big(\frac{t}{s}\Big)^{\lambda(a+c)} \drm s \drm t + 2\frac{c(a+1)}{a+c}\int_0^1\int_0^t s\drm s \drm t \\
& = \frac{2a(1-c^2)(c+1)}{3(2+c-a)(a+c)(1-2a-c)} + \frac{c(a+1)}{3(a+c)}\\
& = \frac{2 - c(c+1 +3ca +3a -2a^2)}{3(2+c-a)(1-2a-c)}.
\end{align*}
\end{proof}
\subsection{The critical regime}
%
\begin{proof}{Theorem}{\ref{CLFT-RERW-CR}}
First, we have from \eqref{QV-N} that for all $t\geq0$
\begin{equation*}
\label{quadNtlog}
\frac{\langle N\rangle_{\floor{n^t}}}{n^t\log n} \longrightarrow 0 \quad\text{a.s.}
\end{equation*}
which implies from Theorem 1.3.24 of \cite{Duflo1997} that
\begin{equation}
\label{Ntlog}
\frac{N_{\floor{n^t}}}{n^t\log n} \longrightarrow 0 \quad\text{a.s.}
\end{equation}
Hereafter, in order to apply Theorem \ref{T-CLTF} to the one-dimensional martingale $(M_n)$, we must once again verify that (H.1), (H.2) and (H.3) are satisfied. \medskip \\
\textbf{(H.1)} Let $w_n = \sqrt{v_n^{-1}}$, we have from \eqref{QV-M}, Remark \ref{remWn} and the fact that $a_{\floor{n^t}}$ is asymtotically equivalent to $n^{-t/2}$ that
\begin{equation*}
w_n\langle M \rangle_{\floor{n^t}}w_n {\underset{n\to\infty}{\longrightarrow}} t(c+1) \quad\text{a.s.}
\end{equation*}
\medskip \\
\textbf{(H.2)}  We also get that Lindeberg's condition is satisfied as $v_n$ is increasing as $\log n$ and we have for all $\eps>0$
\begin{align*}
\lim_{n \rightarrow \infty}  \frac{1}{v_n}\sum_{k=1}^{n^t}\E\big[\Delta M_k^2 \1_{\{|\Delta M_k |>\eps\sqrt{v_n}\}}\big|\cF_{k-1}\big] & \leq \lim_{n\to\infty}\frac{1}{\eps^2v_n^2}\sum_{k=1}^{\floor{n^t}}\E\big[\Delta M_k ^4 \big] \\
&\leq\lim_{n\to\infty}\Big(\frac{v_{\floor{n^t}}}{v_n}\Big)^2\frac{1}{\eps^2 v_{\floor{n^t}}^2}\sum_{k=1}^{\floor{n^t}}  \E\big[\Delta M_k ^4 \big]\\
& \leq\lim_{n\to\infty}\frac{t^2}{\eps^2(\log n^t)^2}\sum_{k=1}^{\floor{n^t}}  \E\big[\Delta M_k ^4 \big].
\end{align*}
Moreover, we have from the very definition of $M_n$ that 
\begin{equation*}
\sum_{k=1}^{n}\E\big[\Delta M_k ^4 \big] = O\Big(\sum_{k=1}^{n} a_k^4\Big) \quad\text{a.s.}
\end{equation*}
and as $a_n$ is asymptotically equivalent to $n^{-1/2}$, we can conclude that
\begin{equation*}
\lim_{n \rightarrow \infty}  \frac{1}{v_n}\sum_{k=1}^{n^t}\E\big[\Delta M_k^2 \1_{\{|\Delta M_k |>\eps\sqrt{v_n}\}}\big|\cF_{k-1}\big] = 0 \quad\text{a.s.}
\end{equation*}
\medskip \\
\textbf{(H.3)} In this particular case, we have $w_t = t (c+1)$. Hence, we obtain that
\begin{equation*}
\big(w_n M_{\floor{n^t}}, \ {t \geq 0}\big) \Longrightarrow \big(W_{t}, \ {t \geq 0}\big)
\end{equation*}
where $W$ is defined as in Theorem \ref{T-CLTF}. 
Moreover, when $a=(1-c)/2$ we obtain from \eqref{DEF-M-N}, \eqref{Ntlog} and the fact that $(a_{\floor{n^t}}v_n)^{-1}$ is asymptotically equivalent to $\sqrt{n^t \log n}^{-1}$ that
\begin{equation*}
\big(\frac{S_{\floor{n^t}}}{\sqrt{n^t\log n}}-\frac{N_{\floor{n^t}}}{\sqrt{n^t\log n}}, \ {t \geq 0}\big) 
\Longrightarrow \frac{1-c}{c+1}\big(W_{t}, \ {t \geq 0}\big). 
\end{equation*}
Consequently, using that $W$ is a centered Brownian motion with variance $(c+1)$, we can conclude  that 
\begin{equation*}
\Big(\frac{S_{\floor{n^t}}}{\sqrt{n^t\log n}}, \ {t \geq 0}\Big) \Longrightarrow 
\sqrt{\frac{(1-c)^2}{c+1}}\big(B_{t}, \ {t \geq 0}\big)
\end{equation*}
and this achieves the proof of Theorem \ref{CLFT-RERW-CR}.
\end{proof}

\section*{Appendix. A non-standard result on martingales}
\addcontentsline{toc}{section}{Appendix. A non-standard result on martingales}
\renewcommand{\thesection}{\Alph{section}}
\renewcommand{\theequation}{\thesection.\arabic{equation}}
\setcounter{section}{1}
\setcounter{equation}{0}
\setcounter{count}{0}

The proofs of our main results rely on the non-standard functional central limit theorem and quadratic strong law for multi-dimensional martingales as for the center of mass of the elephant random walk \cite{LaulinBercu2020}.
A simplified version of Theorem 1 part 2) of Touati \cite{Touati1991} is as follows.
\begin{theorem}
\label{T-CLTF}
Let $(\cM_n)$ be a locally square-integrable martingale of $\R^\delta$ adapted to a filtration $(\cF_n)$,  
with predictable quadratic variation $\langle \cM \rangle_n$.
Let $(V_n)$ be a sequence of non-random square matrices of order $\delta$ such that $\| V_n \|$ 
decreases to $0$ as $n$ goes to infinity. Moreover let $\tau : \R_+ \to \R_+$ be a non-decreasing function going to infinty at infinity. Assume that there exists a symmetric and positive semi-definite matrix $V_t$ that is deterministic and
such that for all $t\geq0$
\begin{equation*}
V_n \langle \cM \rangle_{\tau(nt)} V_n^T \overset{\displaystyle \IP}{\underset{n\to\infty}{\longrightarrow}} V_t. \leqno (\textnormal{H.1})
\end{equation*}
Moreover, assume that Lindeberg's condition is satisfied, that is for all $t\geq 0$ and $\eps >0$,
\begin{equation*}
\sum_{k=1}^{\tau(nt)} \E\bigl[\|V_n \Delta \cM_k \|^2 \1_{\{\|V_n\Delta \cM_k \|>\eps\}}\bigl|\cF_{k-1}\bigr] \overset{\displaystyle \IP}{\underset{n\to\infty}{\longrightarrow}}   0 \leqno (\textnormal{H.2})
\end{equation*}
where $\Delta \cM_n=\cM_n - \cM_{n-1}$. Finally, assume that 
\begin{equation*}
V_t = \sum_{j=1}^q t^{\alpha_j} K_j \leqno (\textnormal{H.3})
\end{equation*}
where $\alpha_j>0$ and $K_j$ is a symmetric matrix, for some $q\in\N^*$.
Then, we have the distributional convergence in the Skorokhod space $D([0,\infty[)$ of right-continuous functions with left-hand limits,
\begin{equation}
\label{CLTF}
\big(V_n \cM_{\tau({nt})},{t \geq 0}\big)  \Longrightarrow  \big(\cW_t, \ {t \geq 0}\big)
\end{equation}
where $\cW=\big(\cW_t, \ {t \geq 0}\big)$ is a continuous $\R^d$-valued centered Gaussian process starting at 0 with covariance, for $0\leq s\leq t$,
\begin{equation}
\label{COV-W-FCLT}
\E[\cW_s\cW_t^T] = V_s.
\end{equation}
\end{theorem}

\vspace{5ex}

\bibliographystyle{acm}
\bibliography{bibliography}

\begin{thebibliography}{10}

\bibitem{baur2019}
{\sc Baur, E.}
\newblock On a class of random walks with reinforced memory.
\newblock {\em J. Stat. Phys.\/} (2020).

\bibitem{Baur2016}
{\sc Baur, E., and Bertoin, J.}
\newblock Elephant random walks and their connection to p{\'o}lya-type urns.
\newblock {\em Physical review. E 94, 052134\/} (2016).

\bibitem{Bercu2018}
{\sc Bercu, B.}
\newblock A martingale approach for the elephant random walk.
\newblock {\em J. Phys. A 51}, 1 (2018), 015201, 16.

\bibitem{BercuHG2019}
{\sc Bercu, B., Chabanol, M.-L., and Ruch, J.-J.}
\newblock Hypergeometric identities arising from the elephant random walk.
\newblock {\em J. Math. Anal. Appl. 480}, 1 (2019), 123360, 12.

\bibitem{BercuLaulin2019}
{\sc Bercu, B., and Laulin, L.}
\newblock On the multi-dimensional elephant random walk.
\newblock {\em J. Stat. Phys. 175}, 6 (2019), 1146--1163.

\bibitem{LaulinBercu2020}
{\sc Bercu, B., and Laulin, L.}
\newblock On the center of mass of the elephant random walk.
\newblock {\em Stochastic Process. Appl. 133\/} (2021), 111 -- 128.

\bibitem{Bertenghi2020}
{\sc Bertenghi, M.}
\newblock Functional limit theorems for the multi-dimensional elephant random
  walk.
\newblock {\em arXiv:2004.02004\/} (2020).

\bibitem{Bertoin2020}
{\sc Bertoin, J.}
\newblock Scaling exponents of step-reinforced random walks.
\newblock {\em Probability Theory and Related Fields 179}, 1 (2021), 295--315.

\bibitem{Businger2018}
{\sc Businger, S.}
\newblock The shark random swim ({L}\'{e}vy flight with memory).
\newblock {\em J. Stat. Phys. 172}, 3 (2018), 701--717.

\bibitem{Chaabane2000}
{\sc Chaabane, F., and Maaouia, F.}
\newblock Th\'{e}or\`emes limites avec poids pour les martingales vectorielles.
\newblock {\em ESAIM Probab. Statist. 4\/} (2000), 137--189.

\bibitem{Coletti2019}
{\sc Coletti, C., and Papageorgiou, I.}
\newblock Asymptotic analysis of the elephant random walk.
\newblock {\em J. Stat. Mech. Theory Exp.}, 1 (2021), 013205.

\bibitem{Coletti2017}
{\sc Coletti, C.~F., Gava, R., and Sch\"{u}tz, G.~M.}
\newblock Central limit theorem and related results for the elephant random
  walk.
\newblock {\em J. Math. Phys. 58}, 5 (2017), 053303, 8.

\bibitem{ColettiN2017}
{\sc Coletti, C.~F., Gava, R., and Sch\"{u}tz, G.~M.}
\newblock A strong invariance principle for the elephant random walk.
\newblock {\em J. Stat. Mech. Theory Exp.}, 12 (2017), 123207, 8.

\bibitem{Duflo1997}
{\sc Duflo, M.}
\newblock {\em Random iterative models}, vol.~34 of {\em Applications of
  Mathematics (New York)}.
\newblock Springer-Verlag, Berlin, 1997.
\newblock Translated from the 1990 French original by Stephen S. Wilson and
  revised by the author.

\bibitem{Fan2020}
{\sc Fan, X., Hu, H., and Xiaohui, M.}
\newblock {Cram{\'e}r moderate deviations for the elephant random walk}.
\newblock {\em J. Stat. Mech. Theory Exp.}, 2 (2021), 023402.

\bibitem{Gonzales2020}
{\sc González-Navarrete, M.}
\newblock Multidimensional walks with random tendency.
\newblock {\em Journal of Statistical Physics 181}, 4 (2020), 1138--1148.

\bibitem{Hall1980}
{\sc Hall, P., and Heyde, C.~C.}
\newblock {\em Martingale limit theory and its application}.
\newblock Academic Press, Inc., New York-London, 1980.
\newblock Probability and Mathematical Statistics.

\bibitem{Janson2004}
{\sc Janson, S.}
\newblock Functional limit theorems for multitype branching processes and
  generalized {P}ólya urns.
\newblock {\em Stochastic Process. Appl. 110}, 2 (2004), 177--245.

\bibitem{Kozam2013}
{\sc Kozma, G.}
\newblock Reinforced random walk.
\newblock In {\em European {C}ongress of {M}athematics}. Eur. Math. Soc.,
  Z\"{u}rich, 2013, pp.~429--443.

\bibitem{Kubota2019}
{\sc Kubota, N., and Takei, M.}
\newblock Gaussian fluctuation for superdiffusive elephant random walks.
\newblock {\em J. Stat. Phys. 177}, 6 (2019), 1157--1171.

\bibitem{Takei2020}
{\sc Miyazaki, T., and Takei, M.}
\newblock Limit theorems for the `laziest' minimal random walk model of
  elephant type.
\newblock {\em J. Stat. Phys. 181}, 2 (2020), 587--602.

\bibitem{Pemantle2007}
{\sc Pemantle, R.}
\newblock A survey of random processes with reinforcement.
\newblock {\em Probab. Surveys 4\/} (2007), 1--79.

\bibitem{Schutz2004}
{\sc Sch\"utz, G.~M., and Trimper, S.}
\newblock Elephants can always remember: Exact long-range memory effects in a
  non-markovian random walk.
\newblock {\em Physical review. E 70, 045101\/} (2004).

\bibitem{Stout1974}
{\sc Stout, W.~F.}
\newblock Maximal inequalities and the law of the iterated logarithm.
\newblock {\em Ann. Probability 1\/} (1973), 322--328.

\bibitem{Touati1991}
{\sc Touati, A.}
\newblock Sur la convergence en loi fonctionnelle de suites de semimartingales
  vers un m\'{e}lange de mouvements browniens.
\newblock {\em Teor. Veroyatnost. i Primenen. 36}, 4 (1991), 744--763.

\bibitem{Vazquez2019}
{\sc V\'azquez~Guevara, V.~H.}
\newblock On the almost sure central limit theorem for the elephant random
  walk.
\newblock {\em J. Phys. A 52}, 1 (2019), 475201.

\end{thebibliography}

\end{document}